\documentclass[12pt]{amsart}
\usepackage{amsmath,amssymb,amsfonts}
\begin{document}
\theoremstyle{plain}
\newtheorem{thm}{Theorem}[section]
\newtheorem{lem}[thm]{Lemma}
\newtheorem{cor}[thm]{Corollary}
\newtheorem{ass}[thm]{Assumption}
\newtheorem{prop}[thm]{Proposition}
\theoremstyle{definition}
\newtheorem{rem}[thm]{Remark}
\newtheorem{claim}[thm]{Claim}
\newtheorem{defn}[thm]{Definition}
\newtheorem{ex}[thm]{Example}
\numberwithin{equation}{section}
\newcommand{\et}{{\rm \acute{e}t}}
\newcommand{\zar}{{\rm zar}}
\newcommand{\dr}{{\rm DR}}
\newcommand{\be}{{\rm B}}
\newcommand{\an}{{\rm an}} 
\newcommand{\red}{{\rm red}}
\newcommand{\codim}{{\rm codim}}
\newcommand{\rank}{{\rm rank\,}}
\newcommand{\Pic}{{\rm Pic}}
\newcommand{\Div}{{\rm Div}}
\newcommand{\Hom}{{\rm Hom}}
\newcommand{\im}{{\rm im}}
\newcommand{\Spec}{{\rm Spec\,}}
\newcommand{\Char}{{\rm char\,}}
\newcommand{\sA}{{\mathcal A}}
\newcommand{\sB}{{\mathcal B}}
\newcommand{\sC}{{\mathcal C}}
\newcommand{\sD}{{\mathcal D}}
\newcommand{\sE}{{\mathcal E}}
\newcommand{\sF}{{\mathcal F}}
\newcommand{\sG}{{\mathcal G}}
\newcommand{\sH}{{\mathcal H}}
\newcommand{\sI}{{\mathcal I}}
\newcommand{\sJ}{{\mathcal J}}
\newcommand{\sK}{{\mathcal K}}
\newcommand{\sL}{{\mathcal L}}
\newcommand{\sM}{{\mathcal M}}
\newcommand{\sN}{{\mathcal N}}
\newcommand{\sO}{{\mathcal O}}
\newcommand{\sP}{{\mathcal P}}
\newcommand{\sQ}{{\mathcal Q}}
\newcommand{\sR}{{\mathcal R}}
\newcommand{\sS}{{\mathcal S}}
\newcommand{\sT}{{\mathcal T}}
\newcommand{\sU}{{\mathcal U}}
\newcommand{\sV}{{\mathcal V}}
\newcommand{\sW}{{\mathcal W}}
\newcommand{\sX}{{\mathcal X}}
\newcommand{\sY}{{\mathcal Y}}
\newcommand{\sZ}{{\mathcal Z}}
\newcommand{\A}{{\mathbb A}}
\newcommand{\B}{{\mathbb B}}
\newcommand{\C}{{\mathbb C}}
\newcommand{\D}{{\mathbb D}}
\newcommand{\E}{{\mathbb E}}
\newcommand{\F}{{\mathbb F}}
\newcommand{\G}{{\mathbb G}}
\renewcommand{\H}{{\mathbb H}}
\newcommand{\J}{{\mathbb J}}
\newcommand{\M}{{\mathbb M}}
\newcommand{\N}{{\mathbb N}}
\renewcommand{\P}{{\mathbb P}}
\newcommand{\Q}{{\mathbb Q}}
\newcommand{\R}{{\mathbb R}}
\newcommand{\T}{{\mathbb T}}
\newcommand{\U}{{\mathbb U}}
\newcommand{\V}{{\mathbb V}}
\newcommand{\W}{{\mathbb W}}
\newcommand{\X}{{\mathbb X}}
\newcommand{\Y}{{\mathbb Y}}
\newcommand{\Z}{{\mathbb Z}}
\catcode`\@=11
\def\opn#1#2{\def#1{\mathop{\kern0pt\fam0#2}\nolimits}} 
\def\bold#1{{\bf #1}}%
\def\underrightarrow{\mathpalette\underrightarrow@}
\def\underrightarrow@#1#2{\vtop{\ialign{$##$\cr
 \hfil#1#2\hfil\cr\noalign{\nointerlineskip}%
 #1{-}\mkern-6mu\cleaders\hbox{$#1\mkern-2mu{-}\mkern-2mu$}\hfill
 \mkern-6mu{\to}\cr}}}
\let\underarrow\underrightarrow
\def\underleftarrow{\mathpalette\underleftarrow@}
\def\underleftarrow@#1#2{\vtop{\ialign{$##$\cr
 \hfil#1#2\hfil\cr\noalign{\nointerlineskip}#1{\leftarrow}\mkern-6mu
 \cleaders\hbox{$#1\mkern-2mu{-}\mkern-2mu$}\hfill
 \mkern-6mu{-}\cr}}}
\let\amp@rs@nd@\relax
\newdimen\ex@
\ex@.2326ex
\newdimen\bigaw@
\newdimen\minaw@
\minaw@16.08739\ex@
\newdimen\minCDaw@
\minCDaw@2.5pc
\newif\ifCD@
\def\minCDarrowwidth#1{\minCDaw@#1}
\newenvironment{CD}{\@CD}{\@endCD}
\def\@CD{\def\A##1A##2A{\llap{$\vcenter{\hbox
 {$\scriptstyle##1$}}$}\Big\uparrow\rlap{$\vcenter{\hbox{%
$\scriptstyle##2$}}$}&&}%
\def\V##1V##2V{\llap{$\vcenter{\hbox
 {$\scriptstyle##1$}}$}\Big\downarrow\rlap{$\vcenter{\hbox{%
$\scriptstyle##2$}}$}&&}%
\def\={&\hskip.5em\mathrel
 {\vbox{\hrule width\minCDaw@\vskip3\ex@\hrule width
 \minCDaw@}}\hskip.5em&}%
\def\verteq{\Big\Vert&&}%
\def\noarr{&&}%
\def\vspace##1{\noalign{\vskip##1\relax}}\relax\let\amp@rs@nd@&\iffalse}\fi
 \CD@true\vcenter\bgroup\relax\let\\=\cr\iffalse}\fi\tabskip\z@skip\baselineskip20\ex@
 \lineskip3\ex@\lineskiplimit3\ex@\halign\bgroup
 &\hfill$\m@th##$\hfill\cr}
\def\@endCD{\cr\egroup\egroup}
\def\>#1>#2>{\amp@rs@nd@\setbox\z@\hbox{$\scriptstyle
 \;{#1}\;\;$}\setbox\@ne\hbox{$\scriptstyle\;{#2}\;\;$}\setbox\tw@
 \hbox{$#2$}\ifCD@
 \global\bigaw@\minCDaw@\else\global\bigaw@\minaw@\fi
 \ifdim\wd\z@>\bigaw@\global\bigaw@\wd\z@\fi
 \ifdim\wd\@ne>\bigaw@\global\bigaw@\wd\@ne\fi
 \ifCD@\hskip.5em\fi
 \ifdim\wd\tw@>\z@
 \mathrel{\mathop{\hbox to\bigaw@{\rightarrowfill}}\limits^{#1}_{#2}}\else
 \mathrel{\mathop{\hbox to\bigaw@{\rightarrowfill}}\limits^{#1}}\fi
 \ifCD@\hskip.5em\fi\amp@rs@nd@}
\def\<#1<#2<{\amp@rs@nd@\setbox\z@\hbox{$\scriptstyle
 \;\;{#1}\;$}\setbox\@ne\hbox{$\scriptstyle\;\;{#2}\;$}\setbox\tw@
 \hbox{$#2$}\ifCD@
 \global\bigaw@\minCDaw@\else\global\bigaw@\minaw@\fi
 \ifdim\wd\z@>\bigaw@\global\bigaw@\wd\z@\fi
 \ifdim\wd\@ne>\bigaw@\global\bigaw@\wd\@ne\fi
 \ifCD@\hskip.5em\fi
 \ifdim\wd\tw@>\z@
 \mathrel{\mathop{\hbox to\bigaw@{\leftarrowfill}}\limits^{#1}_{#2}}\else
 \mathrel{\mathop{\hbox to\bigaw@{\leftarrowfill}}\limits^{#1}}\fi
 \ifCD@\hskip.5em\fi\amp@rs@nd@}
\newenvironment{CDS}{\@CDS}{\@endCDS}
\def\@CDS{\def\A##1A##2A{\llap{$\vcenter{\hbox
 {$\scriptstyle##1$}}$}\Big\uparrow\rlap{$\vcenter{\hbox{%
$\scriptstyle##2$}}$}&}%
\def\V##1V##2V{\llap{$\vcenter{\hbox
 {$\scriptstyle##1$}}$}\Big\downarrow\rlap{$\vcenter{\hbox{%
$\scriptstyle##2$}}$}&}%
\def\={&\hskip.5em\mathrel
 {\vbox{\hrule width\minCDaw@\vskip3\ex@\hrule width
 \minCDaw@}}\hskip.5em&}
\def\verteq{\Big\Vert&}
\def\novarr{&}
\def\noharr{&&}
\def\SE##1E##2E{\slantedarrow(0,18)(4,-3){##1}{##2}&}
\def\SW##1W##2W{\slantedarrow(24,18)(-4,-3){##1}{##2}&}
\def\NE##1E##2E{\slantedarrow(0,0)(4,3){##1}{##2}&}
\def\NW##1W##2W{\slantedarrow(24,0)(-4,3){##1}{##2}&}
\def\slantedarrow(##1)(##2)##3##4{%
\thinlines\unitlength1pt\lower 6.5pt\hbox{\begin{picture}(24,18)%
\put(##1){\vector(##2){24}}%
\put(0,8){$\scriptstyle##3$}%
\put(20,8){$\scriptstyle##4$}%
\end{picture}}}
\def\vspace##1{\noalign{\vskip##1\relax}}\relax\let\amp@rs@nd@&\iffalse}\fi
 \CD@true\vcenter\bgroup\relax\let\\=\cr\iffalse}\fi\tabskip\z@skip\baselineskip20\ex@
 \lineskip3\ex@\lineskiplimit3\ex@\halign\bgroup
 &\hfill$\m@th##$\hfill\cr}
\def\@endCDS{\cr\egroup\egroup}
\newdimen\TriCDarrw@
\newif\ifTriV@
\newenvironment{TriCDV}{\@TriCDV}{\@endTriCD}
\newenvironment{TriCDA}{\@TriCDA}{\@endTriCD}
\def\@TriCDV{\TriV@true\def\TriCDpos@{6}\@TriCD}
\def\@TriCDA{\TriV@false\def\TriCDpos@{10}\@TriCD}
\def\@TriCD#1#2#3#4#5#6{%
\setbox0\hbox{$\ifTriV@#6\else#1\fi$}
\TriCDarrw@=\wd0 \advance\TriCDarrw@ 24pt
\advance\TriCDarrw@ -1em
\def\SE##1E##2E{\slantedarrow(0,18)(2,-3){##1}{##2}&}
\def\SW##1W##2W{\slantedarrow(12,18)(-2,-3){##1}{##2}&}
\def\NE##1E##2E{\slantedarrow(0,0)(2,3){##1}{##2}&}
\def\NW##1W##2W{\slantedarrow(12,0)(-2,3){##1}{##2}&}
\def\slantedarrow(##1)(##2)##3##4{\thinlines\unitlength1pt
\lower 6.5pt\hbox{\begin{picture}(12,18)%
\put(##1){\vector(##2){12}}%
\put(-4,\TriCDpos@){$\scriptstyle##3$}%
\put(12,\TriCDpos@){$\scriptstyle##4$}%
\end{picture}}}
\def\={\mathrel {\vbox{\hrule
   width\TriCDarrw@\vskip3\ex@\hrule width
   \TriCDarrw@}}}
\def\>##1>>{\setbox\z@\hbox{$\scriptstyle
 \;{##1}\;\;$}\global\bigaw@\TriCDarrw@
 \ifdim\wd\z@>\bigaw@\global\bigaw@\wd\z@\fi
 \hskip.5em
 \mathrel{\mathop{\hbox to \TriCDarrw@
{\rightarrowfill}}\limits^{##1}}
 \hskip.5em}
\def\<##1<<{\setbox\z@\hbox{$\scriptstyle
 \;{##1}\;\;$}\global\bigaw@\TriCDarrw@
 \ifdim\wd\z@>\bigaw@\global\bigaw@\wd\z@\fi
 \mathrel{\mathop{\hbox to\bigaw@{\leftarrowfill}}\limits^{##1}}
 }
 \CD@true\vcenter\bgroup\relax\let\\=\cr\iffalse}\fi
 \tabskip\z@skip\baselineskip20\ex@
 \lineskip3\ex@\lineskiplimit3\ex@
 \ifTriV@
 \halign\bgroup
 &\hfill$\m@th##$\hfill\cr
#1&\multispan3\hfill$#2$\hfill&#3\\
&#4&#5\\
&&#6\cr\egroup%
\else
 \halign\bgroup
 &\hfill$\m@th##$\hfill\cr
&&#1\\%
&#2&#3\\
#4&\multispan3\hfill$#5$\hfill&#6\cr\egroup
\fi}
\def\@endTriCD{\egroup}
\title[Germs of de Rham cohomology classes]{Germs of de Rham 
cohomology classes which vanish at the generic point}  
\author{ H\'el\`ene Esnault } 
\address{ Universit\"at Essen, FB6 Mathematik, 45 117 Essen, Germany}
\email{ \begin{tabular}[t]{l} esnault@uni-essen.de\\
viehweg@uni-essen.de \end{tabular}}
\author{ Eckart Viehweg } 
\thanks{ This work has been partly supported by the DFG Forschergruppe
''Arithmetik und Geometrie'' and by the I.H.E.S., Bures-sur-Yvette}
\maketitle
\begin{quote}\footnotesize
{\bf Abstract.} We show that hypergeometric differential
equations, unitary and Gau{\ss}-Manin connections give
rise to de Rham cohomology sheaves which do not admit a 
Bloch-Ogus resolution \cite{BO}. 
The latter is in contrast to Panin's theorem
\cite{P} asserting that corresponding \'etale cohomology sheaves
do fulfill Bloch-Ogus theory.
\end{quote}
\begin{center}
{\bf Germes de classes de cohomologie de de Rham qui
s'annulent au point g\'en\'erique}
\end{center}
\begin{quote}\footnotesize
{\bf R\'esum\'e.} 
Nous montrons que les syst\`emes d'\'equations
hyperg\'eom\'etriques, les connexions unitaires et de
Gau{\ss}-Manin donnent lieu \`a des faisceaux de cohomologie de
de Rham qui n'ont pas de r\'esolution de Bloch-Ogus \cite{BO}. Ce dernier
exemple contraste avec le th\'eor\`eme de Panin \cite{P}
affirmant que des faisceaux semblables en cohomologie
\'etale v\'erifient la th\'eorie de Bloch-Ogus.
\end{quote}
\ \\
{\small {\bf Version abr\'eg\'ee en fran\c{c}ais.} 
Soit $(E, \nabla)$ une connexion plate sur une vari\'et\'e
lisse $S$ sur un corps $k$ alg\'ebriquement clos en
caract\'eristique 0. La restriction des faisceaux de cohomologie
$\sH^i_{DR}((E, \nabla))$ \`a leur valeur au point g\'en\'erique
de $S$ est trivialement injective pour $i=0,1$. Afin de montrer
que pour les exemples de $(E, \nabla)$ \'evoqu\'es plus haut,
cela n'est plus n\'ecessairement le
cas pour $i=2$, nous for\c{c}ons l'existence de germes de sections
de la fa\c{c}on suivante (voir \cite{E}). On remplace $S$ par l'\'eclatement
d'un point. Cela introduit un diviseur exceptionnel
sur lequel la connexion est triviale, et donc acquiert des
sections. Par le morphisme de Gysin, ces sections fournissent
des sections non nulles dans $H^2_{DR}(S, (E, \nabla))$ qui en
particulier s'annulent au point g\'en\'erique de $S$. 
 Que ces sections ne s'annulent pas dans le germe du faisceau en
un point du diviseur exceptionnel provient d'une hypoth\`ese
convenable de r\'esidus dans le cas hyperg\'eom\'etrique, de 
g\'en\'ericit\'e dans le cas unitaire, et de la th\'eorie de
Hodge dans le cas de Gau{\ss}-Manin pour une famille \`a forte 
variation.}\\[.2cm]
\underline{\hspace{4cm}}\\[.4cm]

Let $S$ be a smooth algebraic variety defined over an
algebraically closed field $k$. Bloch-Ogus theory \cite{BO}
provides a acyclic resolution of the Zariski sheaves of
\'etale cohomology $\sH_\et^i(\Z/n(j))$ if $\Char k=0$ or if
$(\Char k,n)=1$, of de Rham cohomology $\sH_\dr^i$ if $\Char
k =0$, and of Betti cohomology $\sH_\be^i$ if $k=\C$. Here
$\sH^i$ denotes the Zariski sheaf associated to the presheaf
$U \mapsto H^i(U)$. The first level of the resolution says that the 
restriction map to the generic point $i_\eta : \eta =\Spec k(S)
\to S$
$$
\sH_\et^i(\Z/n(j)) \>>> i_{\eta *} H_\et^i(\eta, \Z/n(j))
$$
is injective (and similarly for de Rham and Betti cohomology).

Bloch-Ogus theory extends in an obvious way to the sheaves of
cohomology $\sH^i(L)$ with values in a local system of complex
vector spaces $L$ of finite monodromy if $k=\C$, or equivalently
\cite{G} to the de Rham cohomology sheaves
$\sH_\dr^i((E,\nabla))$ of a locally free sheaf $E$ with a flat connection
$\nabla$, the monodromy of which is finite with respect to one
embedding $k \subset \C$ (and hence to all).

A remarkable generalization of the Bloch-Ogus theory for \'etale
cohomology has been given by I. Panin \cite{P}. Let $f:X\to S$
be a projective smooth morphism and let $L$ be a local system of
free $\Z/n$-modules of finite rank (where $n$ is prime to
$\Char k$ if $\Char k >0$). Then the Zariski sheaves
$\sH_\et^i(f,L(j))$ associated to the presheaves
$$
U \longmapsto H_\et^i(f^{-1}(U),L(j))
$$
have a Bloch-Ogus acyclic resolution on $S$. In particular
the restriction to the generic point
$$
\sH_\et^i(f,L(j)) \>>> i_{\eta *} H_\et^i(f^{-1}(\eta),L(j))
$$
is injective, as in the classical case ``$f= {\rm identity}$ and
$L=\Z/n$''. Panin's proof strongly relies on the finiteness of the
local systems involved.

This raises the question of a similar theorem for the de Rham
cohomology in characteristic zero. In this note we give negative
examples:
\subsection{} Bundles $E$ with a flat connection $\nabla$ for which
\begin{equation}\label{inj1}
\sH_\dr^2((E,\nabla)) \>>> i_{\eta *} H_\dr^2(\eta, (E,\nabla))
\end{equation}
is not injective, or equivalently (over $\C$), local systems $L$ of complex
vector spaces for which 
$\sH^2(L) \to i_{\eta *} H^2(\eta,L)$
is not injective (see \ref{ex3} and \ref{ex4}).
\subsection{} Smooth projective morphisms $f:X \to S$ for which 
\begin{equation}\label{inj2}
\sH_\dr^4(f) \>>> i_{\eta *} H_\dr^4(f^{-1}(\eta))
\end{equation}
is not injective (see \ref{ex1} and \ref{ex2}). 
Here $\sH_\dr^i(f)$ denotes the Zariski sheaf
associated to $U \mapsto H_\dr^i(f^{-1}(U)).$
Or equivalently, over $\C$,
$$
\sH_\be^4(f) \>>> i_{\eta *}
H_\be^4(f^{-1}(\eta))=i_{\eta *}\lim_{\>>\hspace*{-.2cm}U\subset 
S\hspace*{-.2cm}>}H_\be^4(f^{-1}(U))$$
is not injective, with a similar notation for Betti cohomology.
\\

Deligne's theorem \cite{D}, saying that $R^\bullet f_* \C =
\oplus_j R^jf_*\C [-j]$ over $k = \C$, together with 
\cite{G} imply that the injectivity of
$$
\sH_\dr^i(f) \>>> i_{\eta *} H_\dr^i(f^{-1}(\eta))
$$
is equivalent to the injectivity 
$$
\sH_\dr^{i-j}((R^jf_*\Omega_{X/S}^\bullet,\nabla)) \>>> i_{\eta *}
H^{i-j}_\dr((\eta,R^jf_*\Omega_{X/S}^\bullet,\nabla))
$$
for all $j$, where $\nabla$ is the Gau{\ss}-Manin connection.
Thus in the second example we will verify that this morphism
is not injective, for $i-j=j=2$.

In all examples given the method to construct a germ of a de
Rham section is as follows.
\subsection{Construction of a section}\label{meth}
Let $(E',\nabla')$ be defined over a smooth variety $S'$ of
dimension at least two, and let $\delta:S \to S'$ be the blow up
of a point $p \in S'$, with exceptional divisor $F$. Let
$(E,\nabla)=\delta^*(E',\nabla')$ be the pullback connection.
Then the restriction 
map 
$$
H_\dr^1(S,(E,\nabla))=H_\dr^1(S',(E',\nabla')) \>>>
H_\dr^1(S'-p,(E',\nabla'))=H_\dr^1(S-F,(E,\nabla))
$$
is an isomorphism. Hence the Gysin map
$$
i_F:H_{\dr,F}^2(S,(E,\nabla))=H_\dr^0(F,(E,\nabla)|_F)=k^{\rank E} \>>> 
H_\dr^2(S,(E,\nabla))
$$
is injective, and any section $i_F(\sigma)$, $\sigma \in
k^{\rank E}$, vanishes at the generic point $\eta$ of $S$.

To show that the maps (\ref{inj1}) and (\ref{inj2}) need not
be injective, we will show that for certain $\sigma$ the image
$i_F(\sigma)$ is non-zero in the stalk
$\sH_\dr^2((E,\nabla))_q$ for all $q\in F$. The latter is equivalent
to saying that for any divisor $C \subset S$, with $F \not\subset C$
\begin{equation}\label{notin}
i_F(\sigma)\not\in i_CH_{\dr,C}^2(S,(E,\nabla)).
\end{equation}
For any smooth dense open subscheme $C_0$ of $C$ and for $S_0=S-(C-C_0)$
one has 
\begin{equation}\label{sm}
H_{\dr,C}^2(S,(E,\nabla))=H_{\dr,C_0}^2(S_0,(E,\nabla))
=H_\dr^0(C_0,(E,\nabla)|_{C_0})=H_\dr^0(\tilde{C},\nu^*(E,\nabla))
\end{equation}
where $\nu:\tilde{C}\to S$ is the normalization of $C \subset S$.
One way to think of this is analytically. Let $L$ be the
kernel of $\nabla$, let $\lambda:S_0 \to S$, $j_0:S_0 - C_0 \to
S_0$ and $j=\lambda\circ j_0$. Then 
\begin{gather*}
H^2_C(S,L)=H^0(S,R^1j_*L) \mbox{ \ \ and}\\
H^2_{C_0}(S_0,L)=H^0(S,\lambda_*R^1j_{0*}L)
=H^0(S_0,R^1j_{0*}L)=H^0(C_0,L|_{C_0}).
\end{gather*}
As $S-S_0$ has codimension $\geq 2$ in $S$, 
$R^i\lambda_*L = 0$ for $i=1,2$.
Thus by the Leray spectral sequence for $j=\lambda\circ j_{0}$
the restriction map $R^1j_{*}L \to
\lambda_*R^1j_{0*}L$ is an isomorphism. Since 
$H^0(C_0,L|_{C_0}) = H^0(\tilde{C},\nu^*L)$
this concludes the proof of (\ref{sm}).
In other words, we do as if $C$ was smooth.

This way to force geometrically the existence of sections which
have nothing to do with the connection was used by the first
author in \cite{E} for example \ref{ex3}. We thank
I. Panin for explaining us the proof of his manuscript \cite{P}.
\section{Gau{\ss}-Manin systems}  
Let $\varphi:Y \to B$ be a semi-stable family of curves of genus $g
\geq 1$ with $B$ a smooth projective curve and $Y$ a smooth
projective surface, defined over an algebraically closed field
$k$ of characteristic zero. We assume throughout this section
\begin{ass}\label{ass}
$\varphi_*\omega_{Y/B}$ is ample.
\end{ass}
Let $B_0$ be the open subscheme of $B$ with $Y_0=\varphi^{-1}(B_0)$
smooth over $B_0$, and as in \ref{meth} let $\delta: S \to
S'=B_0 \times B_0$ be the blow up of a point $p=(b_1,b_2) \in
B_0\times B_0$ with exceptional divisor $F$. Let 
$$
f:X=(Y_0\times Y_0)\times_{(B_0\times B_0)}S \>>> S
$$
be the pullback family. We consider the Gau{\ss}-Manin
connection $(R^2f_* \Omega^\bullet_{X/S},\nabla)$. On the de Rham
cohomology
$$
H^0_\dr(F,(R^2f_* \Omega^\bullet_{X/S},\nabla)|_F) = H^2_\dr(Y_{b_1}\times
Y_{b_2}) 
$$
one has the $F$-filtration $F^0 \supset F^1 \supset F^2$ which
defines a pure Hodge-structure after base extension from
$k$ to $\C$. 
\begin{claim}\label{cl1}
For $\sigma \in \{F^0-F^1\}H^2_\dr(Y_{b_1}\times Y_{b_2})$
and for all $q\in F$ the image $i_F(\sigma)$ is non-zero in the
stalk $\sH_\dr^2((R^2f_*\Omega^\bullet_{X/S},\nabla))_q$. 
\end{claim}
\begin{proof} 
Let $C \subset S$ be a reduced curve with $F \not\subset C$, let
$\nu:\tilde{C}\to S$ be the normalization of $C$, let
$\bar{\delta}: \bar{S} \to B\times B$ be the blow up of $p$,
and let $n:\Gamma \to \bar{S}$ be the normalization of the
closure of $C$ in $\bar{S}$. Let us denote by
\begin{gather*}
\bar{h}:X_\Gamma= (Y\times Y)\times_{(B\times B)} \Gamma \>>>
\Gamma, \\ 
h:X_{\tilde{C}} = (Y_0\times Y_0)\times_{(B_0\times B_0)}
\tilde{C} \>>> \tilde{C},\\
X_F = X\times_S F=\P^1\times Y_{b_1}\times Y_{b_2} \>>> \P^1
\end{gather*}
the induced families of surfaces. The Gysin map
$i_{\tilde{C}}:H^2_\dr(X_{\tilde{C}}) \to 
H^4_\dr(X)$, followed by the restriction map 
$$
\rho_F:H^4_\dr(X) \>>> H^4_\dr(X_F)=H^4_\dr(\P^1\times
Y_{b_1}\times Y_{b_2})
$$
equals the restriction map
$$
\rho_{\nu^{-1}(F)}:H^2_\dr(X_{\tilde{C}}) \>>> \bigoplus_{c\in
\nu^{-1}(F)} H^2_\dr(X_{\tilde{C}}\times_{\tilde{C}}\{c\})=
\bigoplus_{c\in \nu^{-1}(F)}H^2_\dr(\{\nu(c)\}\times
Y_{b_1}\times Y_{b_2}) 
$$
followed by the sum of the Gysin maps
$$
\bigoplus_{c\in \nu^{-1}(F)}H^2_\dr(\{\nu(c)\}\times
Y_{b_1}\times Y_{b_2}) \>>> H^4_\dr(\P^1\times Y_{b_1}\times
Y_{b_2}). 
$$
On the other hand, $\rho_{\nu^{-1}(F)}$ factors through the
surjective map
$$
\rho:H^2_\dr(X_{\tilde{C}}) \>>>
H^0_\dr(\tilde{C},(R^2h_*\Omega^\bullet_{X_{\tilde{C}}/\tilde{C}},
\nabla))=H^2_{\dr,C}(S,(R^2f_*\Omega^\bullet_{X/S},\nabla)),
$$
and $\rho_F \circ i_{\tilde{C}} = \rho'_F \circ i_C \circ \rho$
where $\rho'_F$ is the restriction map
$$
\rho'_F:H^2_\dr(S,(R^2f_*\Omega^\bullet_{X/S},\nabla)) \>>>
H^2_\dr(F,(R^2f_*\Omega^\bullet_{X/S}), \nabla)|_F) =
H^4_\dr(\P^1\times Y_{b_1} \times Y_{b_2}),
$$
and \hfill $
i_C:H^2_{\dr,C}(S,(R^2f_*\Omega^\bullet_{X/S},\nabla))
\>>> H^2_\dr(S,(R^2f_*\Omega^\bullet_{X/S},\nabla))
$ \hspace*{\fill}\\[.3cm]
as in \ref{meth}.
By definition, $\bar{h}$ is a semi-stable family of surfaces,
with singular fibres $Z=\bar{h}^{-1}(\infty)$ for 
$\infty = n^{-1}(\bar{S}-S)=\Gamma-\tilde{C}$. Hence 
$R^2\bar{h}_*\Omega^\bullet_{X_\Gamma /\Gamma}(\log Z)$ is the
Deligne extension \cite{D0} of its restriction to $\tilde{C}$. Thus 
$$
H^0_\dr(\tilde{C},(R^2h_*\Omega_{X_{\tilde{C}}/\tilde{C}}^\bullet,\nabla))=
H^0(\Gamma,\Omega_\Gamma^\bullet(\log \infty)\otimes
R^2\bar{h}_*\Omega_{{X}_\Gamma/\Gamma}^\bullet(\log Z)),
$$
and \hfill
$
\rho(\{F^0/F^1\}H^2_\dr(X_{\tilde{C}}))=\rho(H^2(X_{\tilde{C}},
\sO_{X_{\tilde{C}}})) \subset
H^0(\Gamma,R^2\bar{h}_*\sO_{X_\Gamma}).
$\hspace*{\fill} \\[.3cm]
The sheaf
$$
R^2\bar{h}_* \sO_{{X}_\Gamma} = 
(\bar{\delta} \circ n)^*(pr_1^* R^1\varphi_*\sO_Y
\otimes pr_2^* R^1\varphi_*\sO_Y)
$$
is dual to $(\bar{\delta} \circ n)^* (pr_1^* \varphi_*\omega_{Y/B}
\otimes pr_2^* \varphi_*\omega_{Y/B}).$
Since $\bar{\delta}\circ n$ is finite, the latter
is ample and $H^0(\Gamma,R^2\bar{h}_*\sO_{X_\Gamma}) = 0.$
Thus, since the restriction and Gysin maps are morphisms of
mixed Hodge structures \cite{DII}, 
\begin{equation}\label{eq1}
\im (\rho_F \circ i_{\tilde{C}}) = \im (\rho'_F\circ i_C)
\subset F^3H^4_\dr(\P^1\times Y_{b_1} \times Y_{b_2}).
\end{equation}
On the other hand, $\rho'_F \circ i_F$ is the multiplication by
$(-1)=\deg(\sO_F(F))$. Hence for $\sigma \in
\{F^0-F^1\}H^2_\dr(Y_{b_1}\times Y_{b_2})$ one has
\begin{equation}\label{eq2}
\rho'_F\circ i_F(\sigma) \in \{F^2-F^3\}H^4_\dr(\P^1\times
Y_{b_1} \times Y_{b_2}).
\end{equation}
(\ref{eq1}) and (\ref{eq2}) imply  
$i_F(\sigma)\not\in i_CH_{\dr,C}^2(S,(E,\nabla))$, and, as
explained in \ref{meth}, this proves the claim \ref{cl1}.
\end{proof}
Of course there are lots of families of semi-stable curves $\varphi:Y
\to B$ which satisfy the assumption \ref{ass}. The most
elementary one is:
\begin{ex}\label{ex1}
Choose $\varphi:Y \to B$ to be a semi-stable non-isotrivial family of
elliptic curves. 
\end{ex}
In fact, some power of the sheaf $\varphi_*\omega_{Y/B}$ is the
pullback of an ample invertible sheaf on the moduli space
$\bar{M}_1=\P^1$ of stable elliptic curves.
\begin{ex}\label{ex2}
There exist smooth families of curves $\varphi:Y \to B$ of genus
$g \geq 3$ over a projective curve $B$ with
$\varphi_*\omega_{Y/B}$ ample, and hence there exist smooth
families of surfaces $f:X \to S$ with $S$ projective, for which
the map (\ref{inj2}) is not injective.
\end{ex}
\begin{proof}
Let $M_{g,3}$ and $A_g$ be the moduli spaces of curves of
genus $g$ with level $3$ structure and of $g$-dimensional
principally polarized abelian varieties, respectively. For
$g\geq 3$ the image of $M_{g,3}$ in the Baily-Borel 
compactification of $A_g$ is a projective manifold whose
boundary has codimension larger than or equal to two. Hence
$M_{g,3}$ has a projective compactification with the same
property. Taking hyperplane intersections one obtains a smooth
projective curve $B$ in $M_{g,3}$, and thereby a smooth family
of curves $\varphi:Y \to B$. 

In order to show that for $B$ in general position,
$\varphi_*\omega_{Y/B}$ is ample, we may assume that $k=\C$. 

The monodromy representation of the fundamental group of $M_{g,3}$
is irreducible, and for $B$ in general position the fundamental
group of $B$ maps surjectively to the one of $M_{g,3}$.
Hence the monodromy representation for $B$ is irreducible as well. 

On the other hand, by \cite{K}, 4.10, the sheaf $\varphi_*\omega_{Y/B}$ is
the direct sum of an ample vector bundle and a vector bundle,
flat with respect to the Gau{\ss}-Manin connection. The
irreducibility of the monodromy representation implies that the
latter is trivial.
\end{proof}
\section{Hypergeometric equations}  
As in \ref{meth} let $S'=\P^2-D$, where $D$ is the union of
three lines $H_1$, $H_2$ and $H_3$ in general position.
\begin{ex}\label{ex3}
Choose $a_1, a_2, a_3 \in k$ such that
\begin{enumerate}
\item[i)] the elements  $1, a_i, a_j \in k$ are $\Q$-linearly
independent, for $1 \leq i < j \leq 3$,
\item[ii)] $a_1+a_2+a_3 = 0$;
\end{enumerate}
for example, $a_1=\sqrt{2}$, $a_2=\sqrt{3}$.

Let $\omega\in H^0(\P^2,\Omega^1_{\P^2}(\log D))$ be the unique
form with ${\rm res}_{H_i}\omega = a_i$ for $i=1,2,3$, and 
$(E',\nabla')= (\sO_{S'},d+\omega)$. As in \ref{meth} consider
the blow up $\delta:S \to S'$ with exceptional divisor $F$ and  
the pullback $(E,\nabla)$. We take a section
$0\neq \sigma \in k=H^0(F,(E,\nabla)|_F)$ and regard its image
$i_F(\sigma)$ under the Gysin map. 
\end{ex}
\begin{claim}\label{cl3}
For all reduced curves $C \subset S$ not containing $F$
one has $H^2_{\dr,C}(S, (E,\nabla)) = 0.$
In particular, $0\neq i_F(\sigma) \in \sH_\dr^2((E,\nabla))_q$
for all $q \in F$.
\end{claim}
\begin{proof}
Let $C_0$ be the smooth locus of $C$.
As in (\ref{sm}) 
$$
H^2_{\dr,C}(S, (E,\nabla)) = H^0_\dr(C_0,(E,\nabla)|_{C_0}),
$$
and since $(E,\nabla)$ has rank one, the claim \ref{cl3}
is equivalent to $(E,\nabla)|_{C_0} \neq (\sO_{C_0},d)$.

Let $\bar{\delta}:\bar{S}\to \P^2$ be the blow up of $p$, and
let $n:\Gamma \to \bar{S}$ be the normalization of the closure
of $C$ in $\bar{S}$. For $\infty = n^{-1}\bar{\delta}^{-1}(D)$,
and for some $m_i \in \N$ sufficiently large, one has
$$
H^0_\dr(C_0,(E,\nabla)|_{C_0}) =
H^0(\Gamma,\Omega_\Gamma^\bullet(\log \infty)\otimes
(\bar{\delta}\circ n)^*(\sO_{\P^2}({\textstyle\sum} m_i H_i), d+\omega)).
$$
The residues of $(\bar{\delta}\circ n)^*(\sO_{\P^2}(\sum m_i
H_i), d+\omega)$ along $x \in \Gamma$ are in
$$
\begin{array}{ll}
(\N-\{0\})\cdot(a_i-m_i)& \mbox{ \ \ for \ } \bar{\delta}(n(x)) \in H_i -
\bigcup_{j\neq i} H_j, \\
(\N-\{0\})\cdot(a_i-m_i)+(\N-\{0\})\cdot(a_j-m_j) & \mbox{ \ \
for \ } \bar{\delta}(n(x)) \in H_i\cap H_j.
\end{array}
$$
Indeed, the residue of $(\sO_{\P^2}(\sum m_i H_i),d+\omega)$
along $H_i$ is $a_i-m_i$, and $\Gamma$ lies on a surface $X$
obtained by a sequence of blow ups $X \to S$. As well known,
if $z$ is a smooth point on a variety $Z$ with local coordinates
$x_1, \ldots , x_n$, and $(E,\nabla)$ a differential equation
defined around $z$ by 
$\omega = \sum_{i=1}^n b_i \frac{dx_i}{x_i} + \eta$
for a regular form $\eta$ on $Z$, then the pullback differential
equation on the blow up of $z$ has residue $\sum_{i=1}^n
b_i$ along the exceptional locus.

By the assumption i) in \ref{ex3}, the residues of
$(\bar{\delta}\circ n)^*(\sO_{\P^2}(\sum m_i H_i), d+\omega)$
can not be in $\Q$, and a fortiori $(E,\nabla)|_{C_0}$ can not
be trivial.
\end{proof}\section{Unitary rank one sheaves}
\begin{ex}\label{ex4}
Let $B_1$ and $B_2$ be two non-isogeneous elliptic curves,
defined over $k=\C$, let
$L_i \in \Pic^0(B_i)$ be non-torsion, and let $\nabla_i$ be the
unique unitary connection on $L_i$. Using the notations
introduced in \ref{meth}, we choose $S'=B_1 \times B_2$
and 
$$
(E',\nabla')=(pr_1^*L_1\otimes pr_2^*L_2, pr_1^* \nabla_1
\otimes pr_2^* \nabla_2).
$$
Then for the pullback $(E,\nabla)$ of $(E',\nabla')$ on the blow up
$S$ of a point $p$, the map (\ref{inj1}) will not be injective.
\end{ex}
In fact, for the exceptional divisor $F$ on $S$ one has
$H^0_\dr(F,(E,\nabla)|_F) = \C$
whereas for all reduced curves $C \subset S$ with $F \not\subset
C$ one has:
\begin{claim}\label{cl4}
$H^2_{\dr,C}(S,(E,\nabla)) = 0.$
\end{claim}
\begin{proof}
Let again $C_0$ denote the smooth locus of $C$ and let $n:\Gamma
\to S$ be the normalization. By (\ref{sm})
\begin{gather*}\label{restr}
H_{\dr,C}^2(S,(E,\nabla))=H_\dr^0(C_0,(E,\nabla)|_{C_0})
=H_\dr^0(\Gamma,n^*(E,\nabla))\subset\bigoplus_j
H^0(\Gamma_j,p_{j,1}^*L_1\otimes p_{j,2}^*L_2) 
\end{gather*}
where $\Gamma_j$ are the irreducible components of $\Gamma$,
and where $p_{j,i}$ denotes the restriction of
$pr_i\circ \delta \circ n:\Gamma \to B_i$ to $\Gamma_j$.
If $p_{j,i}$ is dominant, the image $B'_i$ of 
$$
p_{j,i}^*:\Pic^0(B_i) \>>> \Pic^0(\Gamma_j)
$$
is isogeneous to $B_i$ and it is the Zariski closure of the
subgroup generated by $p_{j,i}^*L_i$. Hence if one of the
projections, say $p_{j,1}$, maps $\Gamma_j$ to a point
$p_{j,1}^*L_1\otimes p_{j,2}^*L_2=p_{j,2}^*L_2$ has no global
section. 

If both, $p_{j,1}$ and $p_{j,2}$ are dominant, the two elliptic
curves $B'_1$ and $B'_2$ are not isogeneous, hence $B'_1 \cap
B'_2$ is finite, and 
$H^0(\Gamma_j,p_{j,1}^*L_1\otimes p_{j,2}^*L_2)=0$.
\end{proof}
\bibliographystyle{plain}
\renewcommand\refname{References}
 
\end{document}